\documentclass[final,leqno,onefignum,onetabnum]{siamart171218}
\usepackage{amssymb}

\usepackage{url}
\usepackage{multirow}
\usepackage{algorithmicx}
\usepackage{algpseudocode}
\usepackage[figuresright]{rotating}
\usepackage{amsfonts}
\usepackage{booktabs} 
\usepackage[toc,page,title,titletoc,header]{appendix}
\newtheorem{thm}{Theorem}[section]
\newtheorem{lem}[thm]{Lemma}

\newtheorem{asmp}[thm]{Assumption}

\usepackage{mathtools}
\newcommand\norm[1]{\left\lVert#1\right\rVert}
\author{Rujun Jiang\thanks{School of Data Science, Fudan University, Shanghai, China, rjjiang@fudan.edu.cn}
\and Duan Li\thanks{School of Data Science, City University of Hong Kong, Hong Kong,  dli226@cityu.edu.hk}}
\title{A linear-time algorithm for generalized trust region subproblems}
\begin{document}

\maketitle

\begin{abstract}
In this paper, we provide the first provable linear-time (in term of the number of non-zero entries of the input) algorithm for approximately solving the generalized trust region subproblem (GTRS) of minimizing a quadratic function over a quadratic constraint under some regularity condition.
Our algorithm is motivated  by and extends a recent linear-time algorithm for the trust region subproblem
by  Hazan
and Koren [Math. Program., 2016, 158(1-2): 363-381]. However, due to the non-convexity
and non-compactness of the feasible region, such an extension is nontrivial. Our main contribution is to demonstrate that under some regularity condition, the  optimal solution is in a compact and convex set and  lower and upper bounds of the optimal value  can be computed in linear time. Using these properties, we  develop  a linear-time algorithm for the GTRS.
\end{abstract}
\begin{keywords} the generalized trust region subproblem, semidefinite programming, linear time complexity, approximation algorithms\end{keywords}

\begin{AMS} 90C20,  90C22, 90C26,  68W25\end{AMS}

\pagestyle{myheadings}
\thispagestyle{plain}
\markboth{Linear-time Algorithm for GTRS}{RUJUN JIANG AND DUAN LI}

\section{Introduction}
We consider in this paper the following generalized trust region subproblem (GTRS),
\begin{eqnarray*}
{\rm(GTRS)}~~~~&\min&f(x):= x^TAx+2a^Tx\label{obj}\\
&{\rm s.t.}&h(x):=x^TBx+2b^Tx+d\leq0,\label{con}
\end{eqnarray*}
where $A$ and $B$ are $n\times n$ symmetric matrices which are not necessarily positive semidefinite, $a,b\in\mathbb R^n$ and $d\in\mathbb R$.

When the constraint in (GTRS) is a unit ball, the problem reduces to the classical trust region subproblem (TRS).
The TRS first arose in trust region methods for nonlinear optimization \cite{conn2000trust} and also finds applications in the least square problems \cite{zhang2010derivative} and robust optimization \cite{ben2014hidden}. Various approaches have been derived  to solve the TRS and its variant with additional linear constraints, see
\cite{martinez1994local,more1983computing,ye1992new,rendl1997semidefinite,sturm2003cones,ye2003new,burer2013second,burer2015trust,yang2013two}.
Recently, Hazan and Koren \cite{hazan2016linear} proposed the first linear-time algorithm (with respect to the nonzero entries in the input)  for the TRS, via a linear-time eigenvalue oracle and a linear-time semidefinite programming (SDP) solver based on approximate eigenvector computations \cite{kuczynski1992estimating}. After that,
Wang and Xia \cite{wang2016linear} and Ho-Nguyen
and Kilinc-Karzan \cite{ho2017second} also presented linear-time algorithms to solve the TRS by applying Nesterov's
accelerated gradient descent algorithm to a convex reformation of the TRS, which can also be obtained in linear time.

As a generalization of the TRS, the GTRS  has received a lot of attentions in the literature. The GTRS also admits its own applications such as time of arrival problems \cite{hmam2010quadratic} and subproblems of consensus ADMM in signal processing \cite{huang2016consensus}.
Numerous methods have been developed for solving   the GTRS under various assumptions; see, for example,
\cite{more1993generalizations,stern1995indefinite,ben1996hidden,sturm2003cones,feng2012duality}. 
  Recently, Ben-Tal and den Hertog \cite{ben2014hidden} showed that if the two matrices in the quadratic forms are
  simultaneously diagonalizable (SD) (see \cite{jiang2016simultaneous} for more details about SD conditions), the GTRS can be then transformed into an equivalent second order cone programming (SOCP)
  formulation and thus can be solved efficiently. Salahi and Taati
\cite{salahi2018efficient} also derived an efficient algorithm for solving (GTRS) under the SD condition. Jiang et al. \cite{jiang2018socp} derived an SOCP reformulation for the GTRS when the problem has a finite optimal value and further derived a closed form solution when the SD condition fails.
Pong and Wolkowicz \cite{pong2014generalized} proposed an efficient algorithm based on extreme generalized eigenvalues of a parameterized
matrix pencil for the GTRS, extending the ERW algorithm for the TRS \cite{rendl1997semidefinite}.
 Recently, Adachi and Nakatsukasa \cite{adachi2017eigenvalue} also developed a novel eigenvalue-based algorithm to solve the GTRS.
Jiang and Li \cite{jiang2019novel} proposed a novel convex reformulation for the GTRS and  derived an efficient first order method to solve the reformulation. However,  there is no linear-time algorithm  for the GTRS, while  Hazan and Koren \cite{hazan2016linear} has already proposed their linear-time algorithm for the TRS.
Although it is more general than the TRS, the GTRS still enjoys hidden convexity  the same as the TRS due to the celebrated S-lemma  \cite{yakubovich1971sprocedure,polik2007survey}. There is also evidence that the closely related generalized eigenvalue problem for a positive definite matrix pencil can be solved in linear time \cite{ge2016efficient}. Then a natural question is whether or not there exists a linear-time algorithm for the GTRS. We offer a positive answer to this question in this paper.


In this paper, we derive a linear-time algorithm, Algorithm \ref{alg:TR}, to approximately solve the GTRS with  high probability. The main difficulties in deriving a linear-time algorithm for the GTRS  comes from the non-convex constraint of the GTRS; while this challenging point is not present in the TRS as the constraint in the TRS is convex.
More specifically, the non-convexity of the constraint implies the unboundedness of the feasible region and makes it hard to derive nontrivial initial lower and upper bounds for the GTRS in linear-time.
These difficulties make the direct generalization of the linear-time algorithm for TRS in  \cite{hazan2016linear} inapplicable to
the GTRS.
By addressing these difficulties, we are able to propose a linear-time algorithm for the GTRS based on the work in \cite{hazan2016linear}.  Moreover, our algorithm also inherits the good property of the algorithm in  \cite{hazan2016linear} that avoids the so called hard case\footnote{If the null space of the Hessian matrix of the Lagrangian function, $A+\lambda^* B$, with  $\lambda^{*} $ being the optimal Lagrangian multiplier of problem (P), is orthogonal to $a+\lambda^{*}  b$, we are in the hard case, which can only happen when $\lambda^*$ is an extreme eigenvalue of the matrix pencil $(A,B)$; otherwise we are in the easy case.} by using approximate
eigenvector methods.

The basic idea in our method is to check the feasibility of the following system and then find an  $\epsilon$ optimal solution with a
binary search over $c$,
\begin{equation}
\left.\begin{array}{lll}x^TAx+2a^Tx&\leq&c\\
x^TBx+2b^Tx+d&\le& 0,\end{array}\right.\label{epsilon}
\end{equation}
where $c\in [l,u]$ and $l$ and $u$ are some lower and upper bounds for (GTRS), respectively.
In the TRS, initial $l$ and $u$ can be trivially estimated in linear time as the objective function continuous on the compact feasible region \cite{hazan2016linear}.
However, a linear-time estimation of the lower and upper bounds are nontrivial in the GTRS.
We propose linear-time subroutines that can find a dual feasible solution that in turn helps identify  a lower bound for the primal problem by the weak duality and an upper bound by constructing a feasible solution to the primal problem.
The heart of the binary search is that if system (\ref{epsilon}) is feasible, then $c$ is an upper bound for (GTRS), otherwise system (\ref{epsilon}) is infeasible and $c$ is a lower bound.
In addition, to apply the linear time SDP\ solver in \cite{hazan2016linear}, we introduce a shift $\epsilon$ to system (\ref{epsilon}), i.e.,
\begin{equation}
\left.\begin{array}{lll}x^TAx+2a^Tx&\leq& c-\epsilon\\
x^TBx+2b^Tx+d&\leq &-\frac{\epsilon}{K}.\end{array}\right.\label{deltaep}
\end{equation}
where $K$ is some parameter that can be estimated in linear-time (to be defined in Lemma \ref{shift}). Introducing the parameter $\epsilon$ in the first inequality  shifts  the value of the objective function  with an error $\epsilon$ and introducing the
parameter $\epsilon/ K$  in the second inequality   shifts the value of the objective function at most $\epsilon$ (Lemma \ref{shift}). We then invoke the linear-time SDP solver in \cite{hazan2016linear} that either returns a vector $x$
satisfying system  (\ref{epsilon}), or correctly declares that the direct SDP relaxation of  (\ref{deltaep})
 is infeasible, i.e., a perturbed version of   (\ref{epsilon}) is infeasible. Then via a binary search over $c$, we demonstrate that we can obtain an approximate optimal solution.
However, there are still issues to address when borrowing the SDP solver in \cite{hazan2016linear};
The SDP solver in \cite{hazan2016linear} requires the feasible set of $X$ to be $\{X:X\succeq0,{\rm tr}(X)\le1\}$, where ${\rm tr}(X)$ denotes the trace of matrix $X$, and the direct SDP relaxation of (\ref{deltaep}) lies in an unbounded feasible region.
We will remedy this by showing that the optimal solution of the GTRS must be in  a compact set and further that the optimal solution of the corresponding SDP relaxation should also be in a compact set.

The rest of the paper is organized as follows.
In Section 2, we propose preliminary results that will be used in our algorithm and analysis.
We then illustrate the subroutines to support our main algorithm in Section 3.
We propose our main algorithm and present its analysis in Section 4. Finally, we conclude our paper in Section 5.

\textbf{Notations}
The notations $A\succ0$ and $A\succeq0$ represent that the symmetric matrix $A$ is positive definite and positive semidefinite, respectively.
We use the notation $x(i)$ to denote the $i$th entry of a vector $x$.
We also denote by $v^*$ the optimal value of problem (GTRS).
 Notation $\norm{A}_2$  denotes the operator norm of matrix $A$  and notation $\|a\|$ denotes the Euclidean $l_2$ norm of a vector $a$.
Notations $\lambda_{\max}(A)$ and  $\lambda_{\min}(A)$ denote the largest and smallest eigenvalues of matrix $A$, respectively.
Let $N$ be the total number of nonzero entries of matrices $A$ and $B$ in (GTRS) and assume, without loss of generality, $N\ge n$.
\section{Main results}
In this section, we  review some basic properties of  the GTRS and some oracles and lemmas that will be used later in our algorithms.
\subsection{Preliminary}
Besides the linear-time algorithm for the TRS, Hazan and Koren \cite{hazan2016linear}  also demonstrated that their algorithm can be extended to the GTRS when the quadratic form in the constraint is positive definite. Now let us consider the more general case where $B$ is indefinite.
To avoid some degenerate cases, we do not discuss the case that $B$ is positive semidefinite but singular.

A central tool in this paper is the following linear-time procedures for approximating eigenvectors of sparse matrices, which is based on \cite{kuczynski1992estimating}.

\begin{lem}[Lemmas 3 and 5 in \cite{hazan2016linear}]
Given symmetric matrix $C\in \mathbb R^{n\times n}$ with $\norm{C}_2\leq\rho$ and nonzero entries $N$ and parameters $\epsilon,\delta>0$, there exists an approximate eigenvector oracle, denoted as  $\Call{ApproxmaxEV}{}$,
 that returns a unit vector $x$ with probability of at least $1-\delta$ such that $x^TCx\geq\lambda_{\max}(C)-\epsilon $ and a scalar $\lambda=x^TCx$ in time
\[O(\frac{N\sqrt{\rho}}{\sqrt\epsilon}\log\frac{n}{\delta}).\]
\end{lem}

With the same principle, this oracle can be used to compute an approximate eigenvector to the smallest eigenvalue. In fact, the negativity of the approximate largest eigenvalue for $-C$ gives the smallest approximate eigenvalue of $C$.
\begin{lem}\label{oraclemin}
Given symmetric matrix $C\in \mathbb R^{n\times n}$ with $\norm{C}_2\leq\rho$ and nonzero entries $N$ and parameters $\epsilon,\delta>0$, there exists an approximate eigenvector oracle, denoted as  $\Call{ApproxminEV}{}$,
that returns a unit vector $x$ with probability of at least $1-\delta$ such that $x^TCx\le\lambda_{\min}(C)+\epsilon $ and a scalar $\lambda=x^TCx$ in time
\[O(\frac{N\sqrt{\rho}}{\sqrt\epsilon}\log\frac{n}{\delta}).\]
\end{lem}

We further make the following assumption to ensure the boundedness and existence of an optimal solution.

\begin{asmp}\label{asmp1}
Nontrivial upper bounds $\rho_A$ and $\rho_B$ for $\|A\|_2$ and $\|B\|_2$ are given, i.e., $\|A\|_2\le\rho_A\le 3\|A\|_2$ and $\|B\|_2\le\rho_B\le 3\|B\|_2$. Matrix $B$ has at least one negative eigenvalue with $\lambda_{\min}(B)\le-\xi<0$ for some $\xi>0.$
For the same $\xi,$ there exists a $\mu$ with $\mu\in (0,1]$ such that $\mu A+(1-\mu) B\succeq{\xi}I$.
\end{asmp}

Note that an upper bound  $\rho_A$ (or $\rho _B$) of $\norm{A}_2$ (or $\norm{B}_2$) with $\|A\|_2\le\rho_A\le3\|A\|_2$ (or $\|B\|_2\le\rho_B\le3\|B\|_2$) can be estimated in linear time with high probability, by Theorem 5 in  \cite{hazan2016linear}, and hence the first statement of Assumption \ref{asmp1} can be made without loss of generality.
The second statement of Assumption  \ref{asmp1} is just a numerically stable consideration for the indefiniteness of matrix $B$.
The last statement of   Assumption  \ref{asmp1}  is closely related to the so called regular case, i.e.,  there exists a $\lambda\ge0$ such that $A+\lambda B\succ0$, under which
(GTRS) is bounded below and admits a unique optimal solution \cite{more1993generalizations,pong2014generalized}. Hence our assumption that requires some positive lower bound for the smallest eigenvalue of the positive definite matrix $A+\lambda B$ is reasonable for numerically stable consideration. We assume in the following that
Assumptions   \ref{asmp1} always holds. For ease of notation, we also define $\tilde\xi=\min\{\xi,1\}$ and
\begin{equation}
\label{eq:phi}
\phi=\rho_A+\rho_B+\|a\|+\|b\|+|d|+1.
\end{equation}

When $B$ is indefinite, we can always find a positive constant $K$ in system (\ref{deltaep}) to ensure that the objective value of (GTRS) shifts at most $\epsilon$ if the constraint shifts $\epsilon/K$ when  system  (\ref{deltaep}) is infeasible. To demonstrate this, let us first recall the celebrated  S-lemma \cite{yakubovich1971sprocedure,polik2007survey}.
\begin{lem}
For both $i=1,2$, let $g_{i}(x)=x^TQ_ix+2p_i^Tx+q_i$, where $Q_i$ is an $n\times n$ symmetric matrix, $p_i\in \mathbb R^n$ and $q_i\in\mathbb R$. Assume that  there exists an $\bar x\in\mathbb R^n$ such that $g_2(\bar x)<0$. Then the following two statements are equivalent:\\
{$\rm (S_1)$} There is no $x\in \mathbb R^n$ such that  $g_{1}(x)\leq0$ and $g_{2}(x)< 0$.\\
{$\rm (S_2)$}  There exists a nonnegative multiplier $\lambda\geq0$ such that $ g_{1}(x)+\lambda g_{2}(x)\geq0, ~\forall x\in\mathbb R^n$.
\end{lem}

Note that the assumption of the existence of an $\bar x\in\mathbb R^n$ such that $g_2(\bar x)<0$ automatically holds when $Q_2$ has at least one negative eigenvalue. Using the S-lemma, we have the following results.
\begin{lem}\label{shift}
 Let $K>-\lambda_{\max}(A)/\lambda_{\min}(B)$.
If the system
(\ref{deltaep}) is infeasible for some $\epsilon>0$, then the following system is also infeasible,
\begin{equation}
\left.\begin{array}{lll}x^TAx+2a^Tx&\le &c-2\epsilon,\\
x^TBx+2b^Tx+d&\leq&0.\end{array}\right.\label{2feas}
\end{equation}
\end{lem}
\proof  Since $K$ is an upper bound for $\lambda$ that satisfies $A+\lambda B\succeq0$, we have $\lambda/K<1$. Because $\lambda_{\min}(B)<0$, there exists $\bar x$ such that $\bar x^TB\bar x+2b^T\bar x+d<0$. By the S-lemma, the infeasibility of (\ref{deltaep}) implies
that there exists $\lambda\ge0$ such that
\begin{eqnarray*}
0&\le& x^TAx+2a^Tx-c+\epsilon+\lambda(x^TBx+2b^Tx+d+\epsilon/K)\\
&=&x^TAx+2a^Tx-c+2\epsilon+\lambda(x^TBx+2b^Tx+d)-\varepsilon ,
\end{eqnarray*}
for all $x\in \mathbb R ^n$, where $\varepsilon=(1-\lambda/K)\epsilon>0$.

Due to $K>\lambda\ge0$, for all $x\in R^n$, we have
\begin{eqnarray*}
&&x^TAx+2a^Tx-c+2\epsilon+\lambda(x^TBx+2b^Tx+d-\varepsilon/K)\\
&\ge& x^TAx+2a^Tx-c+2\epsilon+\lambda(x^TBx+2b^Tx+d)-\varepsilon \ge0.
\end{eqnarray*}
Hence for all $\lambda\ge0$, we have
\[
x^TAx+2a^Tx-c+2\epsilon+\lambda(x^TBx+2b^Tx+d-\varepsilon/K)\ge  0,~\forall x.
\]
Moreover, because $\lambda_{\min}(B)<0$, there exists $\bar x$ such that $\bar x^TB\bar x+2b^T\bar x+d-\varepsilon/K<0$. From the S-lemma we conclude that the system
\begin{equation*}
\left.\begin{array}{lll}x^TAx+2a^Tx&\le& c-2\epsilon,\\
x^TBx+2b^Tx+d-\varepsilon/K&<&0\end{array}\right.
\end{equation*}
is infeasible.
This further implies that the system (\ref{2feas}) is infeasible.
\endproof
The above lemma shows that if (\ref{deltaep}) is infeasible,  then $c-2\epsilon$ is a lower bound for (GTRS). Since  an upper bound of $\lambda_{\max}(A)$ is assumed in  Assumption \ref{asmp1}, i.e.,    $\lambda_{\max}(A)\le\rho_A$, and $\lambda_{\min}(B)$ can  be approximately estimated in linear time, parameter $K$ can be estimated in linear time. From Assumption \ref{asmp1}, we have $\mu\lambda_{\max}(A)+(1-\mu)\lambda_{\min}B\succ0$ due to $\mu A+(1-\mu) B\succ0$.
This, together with $\|A\|_2\le \rho_A$ and $\lambda_{\min}(B)\le-\xi$, gives rise to $0<-\lambda_{\max}(A)/\lambda_{\min}(B)\le\rho_A/\xi$ and thus the estimated $K$ can be restricted to be upper bounded by some positive constant (e.g., $\frac{\rho_A}{\xi}+1$).
Hence we set $K=\frac{\rho_A}{\xi}+1$ in our algorithm for simplicity.

\section{Subroutines}
In this section, we present several subroutines to support our main algorithm, Algorithm \ref{alg:TR}.
\subsection{Algorithm to compute  parameter $\mu$ such that $\mu A+(1-\mu) B\succeq\frac{\xi}{2}I$}
In this subsection, a bisection algorithm, Algorithm \ref{alg:pen}, is proposed to find a $\mu$ with $\mu\in(0,1]$ such that $\mu A+(1-\mu) B\succeq\frac{\xi}{2}I$ under Assumption \ref{asmp1}, which is of linear time.
Such a $\mu$ helps us find a lower bound for problem (GTRS) and a compact set in which the optimal solution is located in the following subsections.
We first identify  an interval $(\mu_1,\mu_2] $ where the targeted $\mu$ is located, if exists.  In  Algorithm  \ref{alg:pen}, we  initialize $(\mu_1,\mu_2]$ as $(0,1]$.
In each step, we invoke the eigenvalue oracle $(\lambda,x)=\Call{ApproxminEV}{\mu A+(1-\mu) B,\xi/4,\delta/T}$   to find an approximate smallest eigenvalue of the midpoint of this interval, where $x$ is a unit vector such that $\lambda=x^TAx$ and $T$ is a prescribed maximum iteration number.
When $\lambda<3\xi/4$, we  cut off half of the interval by eliminating either $(\mu_1,\mu]$ if $x^TAx>x^TBx$, or $(\mu,\mu_2]$ if $x^TAx\le x^TBx$. The intuition of this step is that for all  $\nu\in(\mu_1,\mu]$, we have $\lambda_{\min}( \nu A+(1-\nu) B)\le  x^T( \nu A+(1-\nu) B)x\le x^T(\mu A+(1-\mu) B)x=\lambda $ when $x^TAx>x^TBx$. This means that the target $\mu$, if exists, must be in the other half of the interval, i.e.,    $(\mu,\mu_2]$.
 The other situation of this step follows the same argument. We prove  in the following theorem that  under Assumption \ref{asmp1}, such a $\mu$ can be found correctly in linear time with high probability.

\begin{algorithm}[!ht]
  \caption{Compute  parameter $\mu$ such that $\mu A+(1-\mu) B\succeq\frac{\xi}{2}I$}\label{alg:pen}
  \begin{algorithmic}[1]
  \Require symmetric $A,B\in\mathbb R^{n\times n}$ with $\norm{A}_2\le\rho_A$ and $\norm{B}_2\le\rho_{B}$,  $\phi$ in \eqref{eq:phi} and $\xi,\delta>0$
  \Ensure $\mu>0$ such that $\mu A+(1-\mu) B\succeq\frac{\xi}{2}I$ and a unit vector $x$ and $\lambda=x^T(A+\mu B)x$ such that $\lambda_{\min}(\mu A+\mu B)\le \lambda\le\lambda_{\min}(\mu A+\mu B)+\xi/4$; output is correct with probability of at least $1-\delta$
    \Function {PsdPencil}{$A,B,\xi,\phi,\delta$}
    \State initialize $T= \log_2\frac{8\phi}{\xi},~\mu_1=0,~\mu_2=1$
    \For {$i=1:T$}
    \State $\mu=(\mu_1+\mu_2)/2$
    \State invoke $(\lambda,x)\leftarrow\Call{ApproxminEv}{\mu A+(1-\mu) B, \xi/4,\delta/T}$
    \If {$\lambda\geq3\xi/4$}
    \State
    \Return ($\mu$, $\lambda$, $x$)
    \ElsIf {$x^TAx>x^TBx$}
    \State update $\mu_1\leftarrow \mu$ \Comment update $x_1\leftarrow x$ for analysis in Theorem \ref{psdthm}
    \Else
    \State update $\mu_2\leftarrow \mu$ \Comment update $x_2\leftarrow x$ for analysis in Theorem \ref{psdthm}
    \EndIf
    \EndFor
    \State \Return ``Assumption \ref{asmp1} fails."
    \EndFunction
\end{algorithmic}
\end{algorithm}

\begin{thm}\label{psdthm}
Suppose that Assumption \ref{asmp1} holds.  \Call{PsdPencil}{$A,B,\xi,\phi,\delta$} takes at most $\log_2\left(\frac{8\phi}{\xi}\right)$ iterations of $\Call{ApproxminEV}{}$ and returns $(\mu,\lambda,x)$ such that  $\lambda=x^T(\mu A+(1-\mu)B)x\ge 3\xi/4$ and $ \lambda_{\min}(\mu A+(1-\mu)B)\ge\lambda-\xi/4$.
The output is correct with probability of at least $1-\delta$ and  the total runtime is
\[O\left(\frac{N\sqrt\phi}{\sqrt\xi}\log\left(\frac{n}{\delta}\log\frac{1}{\xi}\right)\log\frac{1}{\xi}\right).\]
\end{thm}
\proof
\textit{Runtime:}~~~~
First note that the algorithm invokes $\Call{ApproxminEV}{}$ for at most $T=\log_{2}\frac{8\phi}{\xi}$ iterations.
In  each iteration,  $\Call{ApproxminEV}{}$ is invoked once for matrix $\mu A+(1-\mu) B$,   and the runtime of other main operations (i.e., the matrix vector products $x^TAx$ and $x^TBx$) is dominated by $\Call{ApproxminEV}{}$.
Since $\norm{\mu A+(1-\mu) B}_2\leq \mu\norm{A}_2+(1-\mu)\norm{B}_2\leq\phi$, the runtime for each call of  $\Call{ApproxminEV}{}$ is
$O(\frac{N\sqrt{\phi}}{\sqrt\xi}\log\frac{nT}{\delta})$  from Lemma \ref{oraclemin}.
Hence due to $T=\log_2\frac{8\phi}{\xi}$  each iteration runs in time
$O\left(\frac{N\sqrt{\phi}}{\sqrt\xi}\log\left(\frac{n}{\delta}\log\frac{\phi}{\xi}\right)\right).$
Thus, the total runtime (with at most $T$ iterations) of \Call{PsdPencil}{$A,B,\xi,\phi,\delta$} is
\[O\left(\frac{N\sqrt{\phi}}{\sqrt\xi}\log\left(\frac{n}{\delta}\log\frac{\phi}{\xi}\right)\log\frac{\phi}{\xi}\right).\]

\textit{Correctness:}~~~~
If the algorithm returns $(\mu,\lambda,x)$ for some iteration $i\leq T$, then the returned $\mu$ is the one as required, i.e.,
$\lambda_{\min}({\mu A+ (1-\mu) B})\geq\lambda-\xi/4$. Now it suffices to prove that under Assumption \ref{asmp1}, the algorithm must terminate in some iteration $i\le T$.

Recall that Assumption \ref{asmp1} states that $\exists \mu_0>0$ such that $ \mu_0 A+(1- \mu_0) B\succeq\xi I$.
Since  $\norm{A}_2\leq \rho_A$ and $\norm{B}_2\leq \rho_B$, we have
$A+ B\preceq(\rho_A+\rho_B)I$. Hence for any
$\varrho\in[-\frac{\xi}{4\phi},\frac{\xi}{4\phi}]$, we have
$(\mu_0+\varrho) A+(1- \mu_0-\varrho) B\succeq\xi I-|\varrho|(\rho_A+\rho_B) I\succeq \frac{3\xi}{4}I$ by noting that $\phi\ge \rho_A+\rho_B$.  So for all $\mu'\in[\mu_0-\frac{\xi}{4\phi},\mu_0+\frac{\xi}{4\phi}]\cap(0,1]$,
we have $\mu' A+(1-\mu') B\succeq\frac{3\xi}{4}I$. And the length of the interval  $[\mu_0-\frac{\xi}{4\phi},\mu_0+\frac{\xi}{4\phi}]\cap(0,1]$ is between $\frac{\xi}{4\phi}$ and $\frac{\xi}{2\phi}$.
From the above analysis we know that under Assumption \ref{asmp1}, the interval length of $\{\mu:\mu A+(1-\mu) B\succeq\frac{3\xi}{4}I\}$ is at least $\frac{\xi}{4\phi}$.

If the algorithm does not terminate in the ``for loop", at the end of the  loop, we have  $x_{1}^TC(\mu_1)x_{1}<3\xi/4$ and $x_2^TC(\mu_2)x_2<3\xi/4$ (note that $x_1$ and $x_2$ are defined in the comments in lines 9 and 11 of Algorithm \ref{alg:pen}, respectively), where $C(\mu_i)=\mu_i A+(1-\mu_i)B~(i=1,2)$. Then $x_1^TC(\mu)x_{1}=\mu x_1^TAx_1+(1-\mu)x_1^TBx_1\le \mu_1 x_1^TAx_1+(1-\mu_1)x_1^TBx_1<3\xi/4$ for $\mu\leq \mu_1$ because $x_1^TAx_1>x_1^TBx_1$.
Furthermore, $\lambda_{\min}(C(\mu))\leq x_{1}^T C(\mu)x_{1}<3\xi/4$ for $\mu\leq \mu_1$. Similarly, we have
$\lambda_{\min}(C(\mu))<3\xi/4$ for $\mu\geq \mu_2$.
So if the algorithm  terminates in line 14, we have $\mu_2-\mu_{1}\leq\frac{\xi}{8\phi}$ as the binary search runs for
$\log_2\frac{8\phi}{\xi}$ iterations and the initial length of $\mu_2-\mu_1=1$.  Since for any $\mu$ outside the
interval $[\mu_1,\mu_2]$, we have  $\lambda_{\min}(\mu A+(1-\mu) B)< 3\xi/4$. The length of interval $\{\mu:\lambda_{\min}(\mu A+(1-\mu) B)\geq3\xi/4\}$ is then at most $\mu_2-\mu_1$, which is less than or equal to $\xi/8\phi$. This contradicts Assumption \ref{asmp1}.


The output is correct with probability  of at least $1-\delta$ (from the union bound)  because that each output of \Call{ApproxminEv}{} is correct with probability of at least $1-\delta/T$ and that the total iteration number is at most $T$.
\endproof

\subsection{Lower and upper bounds}
In this subsection, we will show that  Algorithm \ref{alg:bounds} supplies an estimation of initial lower and upper bounds for the optimal value in linear time. The main principle here is that an upper bound can be found by a feasible solution of  (GTRS) and a lower bound can be found by a feasible solution of the Lagrangian dual problem of (GTRS).

\begin{algorithm}[!ht]
\caption{Compute initial lower and upper bounds for (GTRS)}\label{alg:bounds}
\begin{algorithmic}[1]
\Require symmetric $A,B\in \mathbb R^n$ with $\norm{A}_2\le\rho_A$ and $\norm{B}_2\le\rho_B$, $a,b\in \mathbb R^n$, $d\in \mathbb R$, $\delta>0$, $\phi$ in \eqref{eq:phi} and $\xi>0$
\Ensure $\mu_0>0$ such that $\mu_0 A+(1-\mu_0) B\succeq\frac{\xi}{2}I$, a unit vector $x$ and $\lambda=x^T(A+\mu_0 B)x$ such that $\lambda_{\min}(\mu_0 A+\mu_0 B)\le \lambda\le\lambda_{\min}(\mu_0 A+\mu_0 B)+\xi/4$ and lower and upper bounds $l$ and $u$  for (GTRS); output is correct with probability of at least $1-\delta$
\Function {Bounds}{$A,B,\xi,a,b,d,\phi,\delta$}
\State invoke $(\mu_{0},\lambda,x)\leftarrow\Call{PsdPencil}{A,B,\xi$,$\phi$,$\delta/2$}
\State define $\bar\lambda=\lambda-\xi/4$, $\nu_0=(1-\mu_0)/\mu_0$ and $p=a+\nu_0 b$
\State set $l= \frac{d\mu_0}{1-\mu_0}-\frac{\mu_0\norm{p}^2}{\bar\lambda}$ \Comment {initial lower bound}
\State invoke  $(\varpi,y)\leftarrow$\Call{ApproxminEV}{$B,\xi/2,\delta/2$}
\State find a root $\alpha$ to \eqref{eq:alpha}
\State set $u=f(\alpha y)$ \Comment{initial upper bound}
\State     \Return $(\mu_0,\lambda, l,u)$
\EndFunction
\end{algorithmic}
\end{algorithm}

An upper bound for problem $\rm(GTRS)$ is given by $u=f(\alpha y)$, where $\alpha$ is the positive solution of the following quadratic inequality and $y$ is the unit eigenvector returned by  $(\varpi,y)\leftarrow$\Call{ApproxminEV}{$B,\xi/2,\delta/2$},
\begin{equation}\label{eq}
(\alpha y)^T B(\alpha y)+2b^T(\alpha y)+d\leq0.
\end{equation} Under  Assumption \ref{asmp1}, we have $\varpi= y^TBy\leq \lambda_{\min}(B)+\xi/2\leq-\xi/2$.
Consider the quadratic equation
\begin{equation}
\label{eq:alpha}
\varpi\alpha^2+2b^Ty \alpha+d=0.
\end{equation}
The root of the above equation is $\alpha=\frac{-2b^Ty\pm \sqrt{4(b^Ty)^2+4\varpi d}}{-2\varpi}$  with $|\alpha|\le \frac{4\|b\|+\sqrt{4|\varpi d|}}{2|\varpi|}\le\frac{4\|b\|}{\xi}+
\frac{\sqrt{2| d|}}{\sqrt{\xi}} $.
Due to  $y^TBy\le \varpi$, either choice of $\alpha$ is feasible for \eqref{eq}. This yields a feasible solution $\alpha y$ for (GTRS).
Then we have $|f(\alpha  y)|=|\alpha^2y^TAy+2\alpha a^Ty|\le \alpha^2\|A\|_2+2|\alpha|\|a\|\le \left(\frac{4\|b\|}{\xi}+
\frac{\sqrt{2|d|}}{\sqrt{\xi}}\right)^2\rho_A+2\left(\frac{4\|b\|}{\xi}+
\frac{\sqrt{2| d|}}{\sqrt{\xi}}\right)\|a\|=O(\frac{\phi^3}{\tilde \xi^2})$, recalling that $\tilde\xi=\min\{\xi,1\}$.
The vector vector product $b^Ty$ can be done in $O(n)$ and finding a feasible solution of the quadratic equation can be done in $O(1)$. The runtime for computing the approximate smallest eigenvalue is $O(\frac{N\sqrt\phi}{\sqrt{\xi}}\log\frac{n}{\delta})$. The output is correct with probability of at least $1-\delta/2.$ 

Next we illustrate  that a lower bound can be found by a feasible solution of the Lagrangian dual problem.
Note that $(\mu_{0},\lambda,x)\leftarrow\Call{PsdPencil}{A,B,\xi,\phi,\delta/2}$ from line 2 in Algorithm \ref{alg:bounds}.
From Theorem \ref{psdthm} and Algorithm \ref{alg:bounds}, $\lambda$ is a $\xi/4 $ approximate smallest eigenvalue satisfying $\lambda\geq3\xi/4$ under Assumption \ref{asmp1} and the optput is  correct with probability of at least $1-\delta/2$.
Since the  smallest eigenvalue of $\mu_{0} A+(1-\mu_{0})B$ satisfies $\lambda_{\min}(\mu_{0} A+(1-\mu_{0})B)\ge \lambda-\xi/4\ge\xi/2$, we have $\mu_0\rho_A+(1-\mu_0)(-\xi)\ge\xi/2$ and thus $1\ge\mu_0\ge\frac{3\xi}{2(\rho_A+\xi)}$.
The Lagrangian dual problem of (GTRS) is
\begin{eqnarray*}
{\rm (L)}~~~~&&\max_{\nu\ge0}\min\limits_{x} f(x)+\nu  h(x)\\
&=&\max_{\nu\ge0}(x(\nu))^T( A+\nu B)x(\nu)+2(a+\nu b)^Tx(\nu)+d\nu,
\end{eqnarray*}
where $x(\nu)$ is the optimal solution of the inner minimization problem. Letting $\nu_0=(1-\mu_0)/\mu_0$, we have  $P=A+\nu_{0} B\succeq \frac{\xi}{2\mu_{0}} I$. So $\nu_{0} $ is a feasible solution for problem $\rm (L)$.
Define $L(\nu):=\min\limits_{x} f(x)+\nu  h(x)$. Thus from the weak duality, we have the following inequality,
\[ L(\nu_{0})\leq  v^{*}.\]
Letting $p=a+\nu_0 b$, a lower bound of $L(\nu_{0})$ can be found by the following formulation,
\begin{eqnarray*}
L(\nu_{0})&=&\min_x(x+P^{-1}p)^TP(x+P^{-1}p)+d\nu_{0}-p^TP^{-1}p\\
&=& d\nu_{0}-p^TP^{-1}p\\
&\geq& d\nu_{0}-\frac{\norm{p}^2}{\lambda_{\min}(P)}\\
&\geq&d\nu_{0}-\frac{2\mu_0\norm{p}^2}{\xi}.
\end{eqnarray*}
Note that $1\ge\mu_0\ge\frac{3\xi}{2(\rho_A+\xi)}$ implies  $0\le\nu_0\le2\rho_A/3\xi$ and further $\|p\|\le\|a\|+2\rho_A\|b\|/3\xi$. Note also that the main time cost for computing the lower bound is in calling        $\Call{PsdPencil}{}$, which runs in  time, by Theorem \ref{psdthm},
\[O\left(\frac{N\sqrt\phi}{\sqrt\xi}\log\left(\frac{n}{\delta}\log\frac{1}{\xi}\right)\log\frac{1}{\xi}\right).\]

To summarise, in this subsection, we have demonstrated that
\[
u-l\le O(\frac{\phi^3}{\tilde\xi^2})+2\rho_A|d|/3\xi+\frac{\mu_0(\|a\|+\|b\|(1/\mu_0-1))^{2}}{2\xi} , \]
i.e.,
\begin{equation}
\label{eq:bd}
u-l\le O\left(\frac{\phi^3}{\tilde \xi^2}\right), \end{equation}
and the total runtime of \Call{Bounds}{$A,B,\xi,a,b,\phi,\delta$} is
\begin{equation}
\label{eq:bdt}
O\left(\frac{N\sqrt\phi}{\sqrt{\xi}}\log\frac{n}{\delta}+\frac{N\sqrt\phi}{\sqrt\xi}\log\left(\frac{n}{\delta}\log\frac{1}{\xi}\right)\log\frac{1}{\xi}\right)=O\left(\frac{N\sqrt\phi}{\sqrt\xi}\log\left(\frac{n}{\delta}\log\frac{1}{\xi}\right)\log\frac{1}{\xi}\right).
\end{equation}
As we call  \Call{PsdPencil}{$A,B,\xi$,$\phi,\delta/2$} and \Call{ApproxminEV}{$B,\xi/2,\delta/2$} once  each, the output of \Call{Bounds}{$A,B,\xi,a,b,\phi,\delta$} is correct with probability  of at least $1-\delta$ from the union bound.

In summary, we have the following theorem.
\begin{thm}
\label{thm:bound}
Let $\epsilon>0$ and $0<\delta<1$. Under Assumption \ref{asmp1}, with probability of at least $1-\delta$,
the algorithm \Call{Bounds}{$A,B,\xi,a,b,d,\phi,\delta$} computes lower and upper bounds for $\rm (GTRS)$ satsifying
\[
u-l\le O(\frac{\phi^3}{\tilde\xi^2}) ,
\]
and the total runtime is,
\[O\left(\frac{N\sqrt\phi}{\sqrt\xi}\log\left(\frac{n}{\delta}\log\frac{1}{\xi}\right)\log\frac{1}{\xi}\right).\]
\end{thm}

\subsection{Identify feasibility of quadratic systems}
This subsection shows that line 6 of Algorithm \ref{alg:TR} returns an upper bound of the Euclidean norm of the optimal solution and using this bound, the subroutine \Call{Feas}{} can help us to solve the GTRS via bisection.

The subroutine  \Call{Feas}{} utilizes  a linear-time SDP solver, $\Call{RelaxSolve,}{}$ developed in \cite{hazan2016linear}, to approximately solve the following feasibility problem:
\begin{eqnarray}\label{SDPr}
D_i\bullet X\geq\epsilon,~~i=1,2,~~~~X\in \mathcal{K},
\end{eqnarray}
where $D_i$, $i$ = 1, 2, are symmetric matrices with $\norm{D_i}_2\leq1$ and $\mathcal K=\{X:X\succeq0,~{\rm tr(X)}\leq1\}$.
\begin{lem}[Theorem 2 in \cite{hazan2016linear}]\label{l3.2}
Given symmetric matrices $D_1,~D_2\in\mathbb R^{n\times n}$ with $\norm{D_i}_2\leq1$ and $\epsilon,\delta>0$, with probability of at least $1-\delta$, \Call{RelaxSolve}{} outputs a matrix $X\in \mathcal K$ of rank 2 that satisfies $D_i\bullet X\geq \epsilon/2,i=1,2$ or correctly declares that (\ref{SDPr}) is infeasible. The algorithm calls the oracle \Call{ApproxmaxEV}{} at most $O(\log\frac{1}{\epsilon})$ times and can be implemented to run in total time
\[O\left(\frac{N}{\sqrt\epsilon}\log\frac{1}{\epsilon}\log\left(\frac{n}{\delta}\log\frac{1}{\epsilon}\right)\right).\]
\end{lem}

The direct SDP relaxation of  (\ref{deltaep}) is
\begin{equation}\label{SDPr2}
P_i\bullet X\geq\epsilon,~i=1,2,~X_{11}=1,~~X\succeq0,
\end{equation}
where $P_1=\begin{pmatrix}c &-a^T\\
    -a&-A
      \end{pmatrix}$, $P_2=K\begin{pmatrix}-d &-b^{T}\\
    -b&-B
      \end{pmatrix}$.
However, the feasible region for problem (GTRS) may not be contained in $\mathcal{K}$ and thus we cannot directly utilize the SDP solver in \cite{hazan2016linear}. To address this issue, we prove in the next theorem that the optimal solution of (GTRS) must be in a compact set.

\begin{thm}
Let $\mu_0$ be returned by line 3, \Call{Bounds}{$A,B,\xi,a,b,\phi,\delta/2$}, of Algorithm \ref{alg:TR}, $R$ be defined by \eqref{eqn:defR}, i.e., line 6 of Algorithm \ref{alg:TR} and $x^*$ be an optimal solution of problem $\rm(GTRS)$.
Then, under Assumption \ref{asmp1}, it holds that $\norm{x^*}\leq R$ with probability of at least $1-\delta/2$.
Moreover, $R\le O(\phi^{3/2}/\tilde\xi^{3/2})$.
\end{thm}
\proof Let the notations  be the same as in Section 3.2. Hence the output is correct with probability of at least $1-\delta/2$. Also we have   $\lambda_{\min}(\mu_0 A+(1-\mu_0) B)\ge \bar\lambda=\lambda-\xi/4\ge \xi/2>0$. This gives rise to  $\lambda_{\min}(A+\nu_0 B)\ge \bar \lambda/\mu_0 $, where $\nu_0=1/\mu_0-1$. From the optimality of $x^*$, we have $f_{1}(x^*)-u+\nu_0 f_2(x^*)\le 0$, where $u$ is an upper bound of (GTRS) derived from line 3 or 11 of Algorithm \ref{alg:TR}, and thus
\[(x^*)^T(A+\nu_0 B)x^*+2(a+\nu_0 b)^{T}x^{*}-u+\nu_0 d\leq 0.\]
Let $\mathcal{X}=\{x:x^T(A+\nu_0 B)x+2(a+\nu_0 b)^{T}x-u+\nu_0 d\leq 0\}$.
Then $\mathcal{X}\subset\{x:\frac{\bar\lambda}{\mu_0} x^Tx+2(a+\nu_0 b)^{T}x-u+\nu_0 d\leq 0\}$. This  further implies \[\mathcal{X}\subset\left\{x:\norm{x}\le\frac{\mu_0\| a+\nu_0 b\|}{\bar\lambda}+ \sqrt {\norm{\frac{\mu_0}{\bar\lambda}\left(a+\nu_0 b\right)}^2+\frac{\mu_0}{\bar\lambda}(u-\nu_0 d)} \right\}.\]
From the arguments in Section 3.2, we have $3\xi/2(\rho_A+\xi)\le \mu_0\le 1$.
Hence, together with $\bar\lambda\ge\xi/2$, $\mu_0\nu_0=1-\mu_0\le1$ and $u\le O(\phi^3/\tilde\xi^2\})$, it holds that
$R\le O(\phi^{3/2}/\tilde\xi^{3/2})$.
\endproof
By defining
$Y^*=\begin{pmatrix}1\\ x^*\end{pmatrix}\begin{pmatrix}1\\ x^*\end{pmatrix}^T/S$, we have
\[{\rm tr}(Y^*)={\rm tr}(\begin{pmatrix}1\\ x^*\end{pmatrix}^T
\begin{pmatrix}1\\ x^*\end{pmatrix}/S)=\frac{\| x^{*}\|^2+1}{S}\le1,\]
i.e.,  $Y^*\in  \mathcal{K}$, where $S=R^2+1$. This motivates us to solve the following SDP system instead of (\ref{SDPr2}),
\begin{equation}\label{SDPr3}
\frac{1}{\kappa}P_i\bullet Y\geq\frac{\epsilon}{\kappa S},~~~i=1,2,~~Y\in \mathcal{K},
\end{equation}
where $\kappa=\max\{\mu_{A},K\mu_B\}$, $\mu_{A}=\rho_A+2\norm{a}+|c|$, $\mu_B=\rho_B+2\norm{b}+|d|$. These parameters make $\norm{P_i/\kappa}_2\le1,i=1,2$ and the optimal solution, if exists, $Y^*\in\mathcal K$.
Therefore, the SDP feasibility problem (\ref{SDPr3}) can be solved by the linear-time SDP solver, $\Call{RelaxSolve}{}$. Then due to Lemma \ref{l3.2}, with probability of at least $1-\delta$, $\Call{RelaxSolve}{\frac{P_1}{\kappa},\frac{P_2}{\kappa},
 \frac{\epsilon}{\kappa S},\delta}$ either declares that (\ref{SDPr3}) is infeasible, which is further equivalent to the infeasibility of  (\ref{SDPr2}), or returns $Y\in\mathcal K$ such that $\frac{1}{\kappa}P_i\bullet Y\geq\frac{\epsilon}{2\kappa S},i=1,2$, which is further equivalent to that $X=SY$ satisfies
 $P_i\bullet X\geq\epsilon/2,~i=1,2,~X\succeq0$. When $\Call{RelaxSolve}{\frac{P_1}{\kappa},\frac{P_2}{\kappa},
 \frac{\epsilon}{\kappa S},\delta}$ returns $Y\in\mathcal K$ such that $\frac{1}{\kappa}P_i\bullet Y\geq\frac{\epsilon}{2\kappa S}$, we further invoke the  $\Call{SZRotation}{}$ algorithm in \cite{hazan2016linear}, which is a variant of the matrix decomposition procedure in Sturm and Zhang \cite{sturm2003cones}, to find a vector $z$ such that $z^TP_iz\geq\epsilon/2r,i=1,2$ with $r={\rm rank}(Y)\le2$ as there are only two inequalities in \eqref{SDPr3}.

\begin{lem}[\cite{hazan2016linear}]
\label{lem:SZ}
Given a decomposition $X =\sum_{i=1}^r x_i x_i^T$ of a positive semidefinite matrix
$X$ of rank $r$ and an arbitrary matrix $M$ with $M\bullet X\ge a$, \Call{SZRotation}{} outputs a
decomposition $X =\sum_{i=1}^r y_{i} y_i^T$ such that $y_i^TMy_i\ge a/r$ for all $i=1,\ldots,r$. The procedure
runs in time $O(Nr)$, where $N\ge n$ is the number of non-zero entries in $M$.
\end{lem}
\begin{algorithm}[!tbp]
  \caption{Find a feasible solution for (\ref{epsilon}) or declare the infeasibility of (\ref{2feas})}\label{alg:subtr}
  \begin{algorithmic}[1]
  \Require symmetric $A,B\in \mathbb R ^{n\times n}$ with $\norm{A}_2\le \rho_A$
and $\norm{B}_2\le \rho_B$, $a,b,c,d,\epsilon,\delta,\mu_A,\mu_B, K,R>0$ \Ensure find a feasible solution $x$ for (\ref{epsilon}) or declare the infeasibility of (\ref{2feas}); output is correct with probability of at least $1-\delta$
    \Function {FEAS}{$A,B,a,b,c,d,\epsilon,\delta,\mu_A,\mu_B, K,R$}
    \State define $S=(R+1)^2$
    \State define $\kappa=\max\{\mu_{A},K\mu_B\}$
   \State define $(n + 1) \times (n + 1)$ symmetric matrices,
   \begin{equation}\label{Qmatrix}
    Q_1=\frac{1}{\kappa}\begin{pmatrix}c &-a^T\\
    -a&-A
      \end{pmatrix}\text{ and }Q_2=\frac{K}{\kappa}\begin{pmatrix}-d &-b^T\\
    -b&-B
      \end{pmatrix}\end{equation}

    \State invoke $\Call{RelaxSolve}{Q_1,Q_2,\epsilon/(\kappa S),\delta}$     \If {$\Call{RelaxSolve}{}$ returns ``infeasible"}
    \State \Return ``infeasible"
    \Else {~~~~\{\Call{RelaxSolve}{} returns $Y$ such that $ Q_i\bullet Y\geq\frac{\epsilon}{2\kappa S},i=1,2$,\}}
    \State invoke $\Call{SZRotation}{Q_1,SY}$ that return $X=\sum_{i=1}^rz_iz_i^T$ as output
    \If {$r=1$}
    \State
    $ z= z_1$
    \Else \Comment{r=2}
    \State find a vector $z\in\{z_1,z_2\ldots,z_r\}$ for which $z^TQ_2z\geq \epsilon/2r$ and let $\tilde{z}=z(2:n+1)$
    \EndIf  \If {$z(1)\neq0$}
    \State    $x=\tilde{z}/z(1)$
    \Else
    \State set $\alpha=\min\{\frac{\kappa \epsilon}{2r(|2b^T\tilde{z}|+|d|)},\frac{\kappa\epsilon}{2Kr(|2a^T\tilde{z}|+|c|)},1\}$; $x=\tilde{z}/\alpha$
    \EndIf
    \State \Return{$x$}

    \EndIf
    \EndFunction
 \end{algorithmic}
\end{algorithm}
\begin{thm}
\label{thmalg3}
 Given the linear-time SDP\ solver \Call{RelaxSolve}{},  parameters $\epsilon,\delta>0$, and  $Q_1$ and $Q_2$ defined in (\ref{Qmatrix}), Algorithm \ref{alg:subtr},  with probability of at least $1-\delta$,
returns a vector $x\in\mathbb R^n$ satisfying system (\ref{epsilon}),
or correctly declares that (\ref{2feas}) is infeasible.
The  total runtime of Algorithm \ref{alg:subtr}   is
\[O\left(\frac{N\phi^3}{\sqrt{\epsilon\tilde\xi^{5}}}\log\left(\frac{n}{\delta} \log\frac{\phi}{\epsilon\tilde\xi }\right)\log\frac{\phi }{\epsilon\tilde\xi}\right).\]
\end{thm}
\proof
\textit{Runtime:}~~~~The main runtime  is in the two sub-algorithms $\Call{RelaxSolve}{}$ and $\Call{SZRotation}{}$ with  $O\left(\frac{N}{\sqrt{\epsilon'}}\log\left(\frac{n}{\delta} \log\frac{1}{\epsilon'}\right)\log\frac{1}{\epsilon'}\right)$ and $O(Nr)$,
respectively, where $\epsilon'=\epsilon/(\kappa S)$ and $r\leq2$ is the rank of $X$ returned by $\Call{RelaxSolve}{}$. Since $O(Nr)$ is dominated by  $O\left(\frac{N}{\sqrt{\epsilon '}}\log\left(\frac{n}{\delta} \log\frac{1}{\epsilon'}\right)\log\frac{1}{\epsilon'}\right)$,
the total runtime is
\[O\left(\frac{N\sqrt{\kappa S}}{\sqrt{\epsilon}}\log\left(\frac{n}{\delta} \log\frac{\kappa S}{\epsilon}\right)\log\frac{\kappa S }{\epsilon}\right)=O\left(\frac{N\phi^3}{\sqrt{\epsilon\tilde\xi^5}}\log\left(\frac{n}{\delta} \log\frac{\phi}{\epsilon\tilde\xi }\right)\log\frac{\phi }{\epsilon\tilde\xi}\right),\]
where we use $\kappa=\max\{\mu_A,K\mu_B\}\le O(\phi^3/\tilde\xi^2)$ and $S=R^2+1\le O(\phi^3/\tilde\xi^3)$. The bound for $\kappa$ follows that  $\mu_A=\rho_A+2\|a\|+|c|\le2 \phi+|u|+|l|\le O(\phi^3/\tilde\xi^2)$, $K=O(\phi/\xi)$ and $\mu_B=\rho_B+2\|b\|+|d|=O(\phi)$.

\textit{Correctness:}~
If $\Call{RelaxSolve}{}$ returns ``infeasible'',  the SDP relaxation (\ref{SDPr3}) is infeasible and thus (\ref{SDPr2}) is infeasible. This  implies the infeasibility of (\ref{deltaep})  since (\ref{SDPr2}) is a relaxation of (\ref{deltaep}). And the infeasibility of (\ref{deltaep}) further implies the infeasibility of (\ref{2feas}) by Lemma \ref{shift}.

Now let us assume that $\Call{RelaxSolve}{}$ returns $Y\in\mathcal K$ such that $Q_i\bullet Y\geq\frac{\epsilon}{2\kappa S},~i=1,2$. Then $X=SY$ satisfies
 $P_i\bullet X\geq\epsilon/2\kappa,~i=1,2,~X\succeq0$. As shown in  $\Call{RelaxSolve}{}$ in \cite{hazan2016linear},
we have $Y=qy_1y_1^T+(1-q)y_2y_2^T$ for $q\in[0,1]$, which means that ${\rm rank(} Y)\le2$ and thus ${\rm rank }(X)\le2$. It follows from Lemma \ref{lem:SZ} that invoking $\Call{SZRotation}{Q_1,X}$ yields a
solution $X=\sum_{i=1}^rz_iz_i^T$  such that  $z_i^TQ_1z_i\geq\epsilon/2r\kappa, ~i=1,2$. Together with  $Q_2\bullet X\geq\epsilon/2\kappa$, we conclude that at least one of
$z_i$  satisfies $z_i^TQ_2z_i\geq\epsilon/2r\kappa $. So $\Call{SZRotation}{Q_1,X}$ indeed finds a vector $z\in\{z_1,z_2\}$ such that $z^TQ_jz\geq\epsilon/2r\kappa, ~(j=1,2)$, which further implies that
\begin{equation*}
\left.\begin{array}{lll}\tilde{z}^TA\tilde{z}+2z(1) a^T\tilde{z}&\le& cz(1)^2 -\epsilon/2r\kappa,\\
\tilde{z}^TB\tilde{z}+2z(1) b^T\tilde{z}+z(1)^2 d&\le&-\epsilon/2rK\kappa,\end{array}\right.
\end{equation*}
where $\tilde{z}=z(2:n+1)$.

If $z(1)\neq0$, by dividing $z(1)^2$ on both sides of the two inequalities in the above system and letting $x=\tilde{z}/z(1)$,
we have
\begin{equation*}
\left.\begin{array}{lll}x^TAx+2a^Tx-c&\le&-\epsilon/2r\kappa z(1)^{2},\\
x^TBx+2b^Tx+d&\le& -\epsilon/2rK\kappa z(1)^{2}.\end{array}\right.
\end{equation*}
Then we have
\begin{equation*}
\left.\begin{array}{lll}x^TAx+2a^Tx&\le& c,\\
x^TBx+2b^Tx+d&\le& 0 \end{array}\right.
\end{equation*}
as required. Else we have $ z(1)=0$. Note that
\begin{eqnarray*}
&&\tilde{z}^TA\tilde{z}\leq -\epsilon/2r\kappa,\\
&&\tilde{z}^TB\tilde{z}\leq -\epsilon/2rK\kappa.
\end{eqnarray*}
By setting $\alpha=\min\{\frac{\epsilon}{2r\kappa(|2a^T\tilde{z}|+|c|)},\frac{\epsilon}{2rK\kappa(|2b^T\tilde{z}|+|d|)},1\}\leq1$ and $x=\tilde{z}/\alpha$,
we have
\begin{eqnarray*}
&&x^TAx=\frac{\tilde{z}^TA\tilde{z}}{\alpha^2}\\
&\leq&-\frac{\epsilon}{2r\kappa\alpha^2}\\
&\leq&-\frac{|2a^T\tilde{z}|+|c|}{\alpha}~~~({\rm due~to}~ -1/\alpha\leq -\frac{2r\kappa(|2a^T\tilde{z}|+|c|)}{\epsilon}) \\
 &\leq&-|2 a^T\tilde{z}/\alpha|-|c|.~~~~~({\rm due~to}~ \alpha\leq1)
\end{eqnarray*} Thus, $x^TAx+2a^Tx-c\leq x^TAx+2| a^T\tilde{z}/\alpha|+|c|\leq 0.$
Similarly, we have
\[x^TBx+2b^Tx+d\leq 0.\]
Hence $x$ is indeed a solution to system (\ref{epsilon}).
\endproof
\section{Main algorithm}
In our main algorithm, Algorithm \ref{alg:TR}, we first invoke a subroutine  \Call{Bounds}{} (defined in Section 3.2) to compute an initial estimation for  lower and upper bounds, $l$ and $u$, for problem (GTRS) and then use  bisection techniques, by invoking \Call{Feas}{} (defined in Section 3.3) for at most
$O\left(\log(\frac{u-l}{\epsilon})\right)$ iterations, to obtain a feasible solution $\tilde{x}$ with $f(\tilde{x})\leq v^*+\epsilon$, where  $v^*$ is an optimal value of (GTRS).

\begin{algorithm}[!ht]
\caption{Find an $\epsilon$ optimal solution for (GTRS)}\label{alg:TR}
\begin{algorithmic}[1]
\Require symmetric $A,B\in \mathbb R^n$ with $\norm{A}_2\le\rho_A$ and $\norm{B}_2\le\rho_B$, $a,b\in \mathbb R^n$, $d\in \mathbb R$, $\epsilon,\delta>0$ and $\xi$  \Ensure an $\epsilon$ optimal solution; output is correct with probability of at least $1-\delta$
\Function {GTRS}{$A, B,a,b,d,\epsilon,\rho_A,\rho_B,\xi,\delta$}
\State let $K=\frac{\rho_A}{\xi}+1$ and $\phi=\rho_A+\rho_B+\|a\|+\|b\|+|d|+1$
\State invoke $(\mu_0,\lambda, l,u)\leftarrow$\Call{Bounds}{$A,B,\xi,a,b,\phi,\delta/2$}~
\State define $\nu_0=1/\mu_0-1$, $\bar\lambda=\lambda-\xi/4$, $c=\frac{l+u}{2}$, $\epsilon'=\frac{\epsilon}{7}$, $T=\log_2(\frac{u-l}{\epsilon'})$, $\delta'=\frac{\delta}{2T}$, $\mu_A=\rho_A+2\norm{a}+|c|$ and $\mu_B=\rho_B+2\norm{b}+|d|$

\For {$t=1:T$}
\State set\begin{equation}\label{eqn:defR}
R=\frac{\mu_{0}\| a+\nu_0 b\|}{\bar\lambda}+ \sqrt {\left(\norm{\frac{\mu_{0}}{\bar\lambda}\left(a+\nu_0 b\right)}^2+\frac{\mu_{0}}{\bar\lambda}(u-\nu_0 d)\right)}
\end{equation}
\State invoke \Call{Feas}{$A,B,a,b,c,d,\epsilon',\xi,\delta',\mu_A,\mu_B, K,R$}
\If {\Call{Feas}{$A,B,a,b,c,d,\epsilon',\xi,\delta',\mu_A,\mu_B,K,R$} returns ``infeasible"}
\State set $l=c-2\epsilon'$; $c=(l+u)/2$
\Else {~\Call{Feas}{$A,B,a,b,c,d,\epsilon',\xi,\delta',\mu_A,\mu_B,K,R$} returns  a feasible solution  $x$ to (\ref{epsilon})}
\State set $u=\min\{u,x^TAx+2a^Tx\}$, $c=(l+u)/2$ and $\mu_A=\rho_A+2\norm{a}+|c|$
\EndIf
\EndFor
\State    \Return{$x$}
\EndFunction
 \end{algorithmic}
\end{algorithm}

We are now ready to present our main result, which shows us the correctness and linear runtime of Algorithm \ref{alg:TR}.
\begin{thm}\label{thm:main}
 Let $\epsilon>0$ and $0<\delta<1$. Under Assumption \ref{asmp1}, with   probability of at least $1-\delta$,
 Algorithm \ref{alg:TR} returns an $\epsilon$ optimal solution $\tilde{x}$ to $(\rm GTRS)$, i.e.,
 a feasible solution $\tilde{x}$ with $f(\tilde{x})\leq v^*+\epsilon$,  where  $v^*$ is an optimal value of $\rm(GTRS)$. The total runtime is
\[O\left(\frac{N\phi^3}{\sqrt{\epsilon\tilde\xi^5}}\log\left(\frac{n}{\delta} \log\frac{\phi}{\epsilon\tilde\xi }\right)\log\frac{\phi }{\epsilon\tilde\xi}\log\frac{\phi}{\epsilon\tilde\xi}\right).\]
 \end{thm}
\proof
\textit{Correctness:}~~~~Section 3.2 shows that
the lower and upper bounds $l$ and $u$ can be estimated, with probability of at least  $1-\delta/2$, by \Call{Bounds}{$A,B,\xi,a,b,\phi,\delta/2$} and $u-l\le O(\phi^3/\tilde\xi^2)$ from \eqref{eq:bd}.

From Theorem \ref{thmalg3}, we know that the subroutine \Call{Feas}{} either returns a feasible solution for (\ref{obj}) (yielding a new  upper bound $u=\min\{c,x^TAx+2a^Tx\}$) or declares the infeasibility of (\ref{deltaep}) (yielding a new lower bound $l=c-2\epsilon'$ by Lemma \ref{shift}).
Now consider the loop in lines 5--13 of Algorithm \ref{alg:TR}. Let $l_p$ and $u_p$ denote the values of $l$ and $u$ in the end of the $p$th iteration in the ``for" loop (particularly, let $l_0$ and $u_0$ be the initial values of $l$ and $u$) and then the length of $u-l$ is at most $\frac{u_p-l_p}{2}+2\epsilon'$ at the end of the current iteration.
At the end of the main loop of Algorithm \ref{alg:TR},  the length of $u_{T}-l_{T}$ satisfies
\[u_T-l_T\leq\frac{u_{0}-l_{0}}{2^T}+(2+1+\cdots+\frac{1}{2^{T-1}})\epsilon'\leq\epsilon'+4\epsilon'=5\epsilon'.\]
From Lemma \ref{shift}, we have $l_T-2\epsilon'\leq f(x^{*})\leq f(x)= u_T$.
Thus $u_T-f(x^{*})\leq u_T-(l_T-2\epsilon')\leq(5+2)\epsilon'=7\epsilon'=\epsilon$. So $f(x)\leq f(x^{*})+\epsilon$.
That is, after
\[O\left(\log\frac{u-l}{\epsilon}\right)=O\left(\log\frac{\phi}{\tilde\xi\epsilon}\right)\] iterations of binary search, we obtain an $\epsilon$ optimal solution.

\textit{Runtime:}~~~~The main runtime of Algorithm \ref{alg:TR} is in subroutines  \Call{Bounds}{} and \Call{Feas}{}. The equation \eqref{eq:bdt} in Section 3.2   shows that  the main operations in  \Call{Bounds}{$A,B,\xi,a,b,\phi,\delta/2$} run in time
\begin{equation}
\label{eq:time1}
O\left(\frac{N\sqrt\phi}{\sqrt\xi}\log\left(\frac{n}{\delta}\log\frac{1}{\xi}\right)\log\frac{1}{\xi}\right).
\end{equation}
The result returned by  \Call{Bounds}{$A,B,\xi,a,b,\phi,\delta/2$} is correct with probability of at least $1-\delta/2$.

Note that Algorithm \ref{alg:TR} invokes \Call{Feas}{}
 $O(\log\frac{u-l}{\epsilon'})$ times. Then from Theorem \ref{thmalg3} in Subsection 3.3,  the total time of lines 5-14 is
\[O\left(\frac{N\phi^3}{\sqrt{\epsilon'\tilde\xi^{5}}}\log\left(\frac{n}{\delta'} \log\frac{\phi}{\epsilon'\tilde\xi }\right)\log\frac{\phi }{\epsilon'\tilde\xi}\log\frac{u-l}{\epsilon'}\right),\]
which is equivalent to,
\begin{equation}
\label{eq:time2}O\left(\frac{N\phi^3}{\sqrt{\epsilon\tilde\xi^{5}}}\log\left(\frac{n}{\delta} \log\frac{\phi}{\epsilon\tilde\xi }\right)\log\frac{\phi }{\epsilon\tilde\xi}\log\frac{\phi}{\epsilon\tilde\xi}\right),
\end{equation}
 by noting   $u-l\le O(\phi^3/\tilde\xi^2)$, $\delta'=\delta/2T$, $T=\log_2((u-l)/\epsilon')$  and $\epsilon'=\epsilon/7$. The output is correct with probability of, by noting that $\delta'=\frac{\delta}{2T}$, at least $1-T\times \frac{\delta}{2T}=1-\frac{\delta}{2}$.

Hence the output of the whole algorithm is correct with probability of at least $1-\delta$.
 Combining the runtime of \eqref{eq:time1} and \eqref{eq:time2}, we conclude that the total runtime is
\[O\left(\frac{N\phi^3}{\sqrt{\epsilon\tilde\xi^{5}}}\log\left(\frac{n}{\delta} \log\frac{\phi}{\epsilon\tilde\xi }\right)\log\frac{\phi }{\epsilon\tilde\xi}\log\frac{\phi}{\epsilon\tilde\xi}\right).\]

\endproof
It is interesting to compare our linear-time algorithm to existing algorithms in the literature.
The algorithms   in \cite{adachi2017eigenvalue,jiang2019novel} and our algorithm   all require  the regularity condition $A+\lambda B\succ0$. However, both   \cite{adachi2017eigenvalue}  and  \cite{jiang2019novel} do not give a way to compute such a $\lambda$ in linear time; While our subroutine \Call{PsdPencil}{} gives a linear-time algorithm for such computation.
Another drawback of the algorithms in   \cite{adachi2017eigenvalue}  and  \cite{jiang2019novel} is that they both require  exact computation of an extreme eigenpair and it is unknown how much the solutions of (GTRS) will be perturbed if an inexact eigenpair was used.
Though in practice,  \cite{jiang2019novel} and  \cite{adachi2017eigenvalue}   use an inexact eigenpair, their theoretical guarantee is missing.
Besides, the algorithm in \cite{jiang2019novel}  requires $O(\frac{L\|x_0-x^*\|^2}{\epsilon})$ iterations to  solve a minimax reformulation of GTRS to achieve an $\epsilon$ optimal solution, where $L$ is the Lipschitz constant for the gradients of functions in the minimax reformulation.
Each iteration needs several matrix vector products and thus  the complexity in each iteration is $O(N)$. So the total complexity is then $O(\frac{L\|x_0-x^*\|^2}{\epsilon})$ (this may be roughly considered as $ O(\frac{\phi^4 }{\epsilon\tilde\xi^3})$ due to $L\le \phi$ and $\|x_0-x^*\|^2\le (2R)^2\le O(\phi^3/\xi^3)$, which is worse than our complexity that is proportional to $\frac{N\phi^3}{\sqrt{\epsilon\tilde \xi^5}}$ as shown in Theorem \ref{thm:main}.
Next let us give comparisons with the results in \cite{ben2014hidden,jiang2018socp,pong2014generalized}. The method in  \cite{ben2014hidden} involves a process in simultaneously diagonalizing two matrices, whose computation is not given there.
The paper \cite{jiang2018socp} considers a case that two matrices are not simultaneously diagonalizable, a case that is numerically unstable, which is not considered in our setting.
The algorithm in \cite{pong2014generalized} solves an extreme generalized eigenpair of a parameterized matrix pencil for the GTRS at each step, whose iteration complexity is unknown. And how an inexact computation of reformulation or extreme eigenpairs will influence the final result  is not provided in   \cite{ben2014hidden,jiang2018socp,pong2014generalized}.
In summary, our methods represent the first provable linear-time algorithm for the GTRS.

\section{Conclusion}
In this paper, we have presented the first linear-time algorithm to approximately solve the generalize trust region subproblem, which extends the recent result in \cite{hazan2016linear} for the trust region subproblem. Our algorithm avoids diagonalization or factorization of matrices as that in \cite{hazan2016linear}.
Our algorithm also has the same time complexity as in the linear-time algorithm for TRS in \cite{hazan2016linear} as well as  in generalized
eigenvector computation \cite{ge2016efficient}.
Similar to \cite{hazan2016linear}, our algorithm avoids the ``hard case" by using an approximate  linear-time SDP solver.
A byproduct of this paper is to provide a linear-time algorithm, \Call{PsdPencil}{}, to detect a $\lambda$ such that $A+\lambda B\succ0$ under mild conditions, which may be of some independent interest for readers.
Our future research will focus on extending the current algorithm to some variants of the GTRS with additional linear constraints or an additional unit ball constraint.

\section*{Acknowledgements}
This research was partially supported by Shanghai Sailing Program 18YF1401700, Natural Science Foundation of China (NSFC) 11801087 and Hong Kong Research Grants Council under
Grants 14213716 and 14202017.  The authors would like to thank the two anonymous referees for the invaluable comments that improve the quality of the paper significantly.
\bibliographystyle{siamplain}
\bibliography{reference}

\begin{thebibliography}{10}

\bibitem{adachi2017eigenvalue}
{\sc S.~Adachi and Y.~Nakatsukasa}, {\em Eigenvalue-based algorithm and
  analysis for nonconvex qcqp with one constraint}, Mathematical Programming,
  (2017), pp.~1--38.

\bibitem{ben2014hidden}
{\sc A.~Ben-Tal and D.~{den Hertog}}, {\em Hidden conic quadratic
  representation of some nonconvex quadratic optimization problems},
  Mathematical Programming, 143 (2014), pp.~1--29.

\bibitem{ben1996hidden}
{\sc A.~Ben-Tal and M.~Teboulle}, {\em Hidden convexity in some nonconvex
  quadratically constrained quadratic programming}, Mathematical Programming,
  72 (1996), pp.~51--63.

\bibitem{burer2013second}
{\sc S.~Burer and K.~M. Anstreicher}, {\em Second-order-cone constraints for
  extended trust-region subproblems}, SIAM Journal on Optimization, 23 (2013),
  pp.~432--451.

\bibitem{burer2015trust}
{\sc S.~Burer and B.~Yang}, {\em The trust region subproblem with
  non-intersecting linear constraints}, Mathematical Programming, 149 (2015),
  pp.~253--264.

\bibitem{conn2000trust}
{\sc A.~R. Conn, N.~I. Gould, and P.~L. Toint}, {\em Trust Region Methods},
  vol.~1, {Society for Industrial and Applied Mathematics (SIAM),
  Philadelphia}, 2000.

\bibitem{feng2012duality}
{\sc J.-M. Feng, G.-X. Lin, R.-L. Sheu, and Y.~Xia}, {\em Duality and solutions
  for quadratic programming over single non-homogeneous quadratic constraint},
  Journal of Global Optimization, 54 (2012), pp.~275--293.

\bibitem{ge2016efficient}
{\sc R.~Ge, C.~Jin, P.~Netrapalli, A.~Sidford, et~al.}, {\em Efficient
  algorithms for large-scale generalized eigenvector computation and canonical
  correlation analysis}, in International Conference on Machine Learning, 2016,
  pp.~2741--2750.

\bibitem{hazan2016linear}
{\sc E.~Hazan and T.~Koren}, {\em A linear-time algorithm for trust region
  problems}, Mathematical Programming, 158 (2016), pp.~363--381.

\bibitem{hmam2010quadratic}
{\sc H.~Hmam}, {\em Quadratic {{Optimization}} with {{One Quadratic Equality
  Constraint}}.}, tech. report, Warfare and Radar Division DSTO Defence Science
  and Technology Organisation, Australia, Report DSTO-TR-2416, 2010.

\bibitem{ho2017second}
{\sc N.~Ho-Nguyen and F.~Kilinc-Karzan}, {\em A second-order cone based
  approach for solving the trust-region subproblem and its variants}, SIAM
  Journal on Optimization, 27 (2017), pp.~1485--1512.

\bibitem{huang2016consensus}
{\sc K.~Huang and N.~D. Sidiropoulos}, {\em Consensus-admm for general
  quadratically constrained quadratic programming}, IEEE Transactions on Signal
  Processing, 64, pp.~5297--5310.

\bibitem{jiang2016simultaneous}
{\sc R.~Jiang and D.~Li}, {\em Simultaneous {{Diagonalization}} of {{Matrices}}
  and {{Its Applications}} in {{Quadratically Constrained Quadratic
  Programming}}}, SIAM Journal on Optimization, 26 (2016), pp.~1649--1668.

\bibitem{jiang2019novel}
{\sc R.~Jiang and D.~Li}, {\em Novel reformulations and efficient algorithm for
  the generalized trust region subproblem}, accepted by SIAM Journal on
  Optimization,  (2019).

\bibitem{jiang2018socp}
{\sc R.~Jiang, D.~Li, and B.~Wu}, {\em {{SOCP}} reformulation for the
  generalized trust region subproblem via a canonical form of two symmetric
  matrices}, Mathematical Programming, 169 (2018), pp.~531--563.

\bibitem{kuczynski1992estimating}
{\sc J.~Kuczy{\'n}ski and H.~Wo{\'z}niakowski}, {\em Estimating the largest
  eigenvalue by the power and lanczos algorithms with a random start}, SIAM
  Journal on Matrix Analysis and Applications, 13 (1992), pp.~1094--1122.

\bibitem{martinez1994local}
{\sc J.~M. Mart{\'\i}nez}, {\em Local minimizers of quadratic functions on
  {{Euclidean}} balls and spheres}, SIAM Journal on Optimization, 4 (1994),
  pp.~159--176.

\bibitem{more1993generalizations}
{\sc J.~J. Mor{\'e}}, {\em Generalizations of the trust region problem},
  Optimization Methods and Software, 2 (1993), pp.~189--209.

\bibitem{more1983computing}
{\sc J.~J. Mor{\'e} and D.~C. Sorensen}, {\em Computing a trust region step},
  SIAM Journal on Scientific and Statistical Computing, 4 (1983), pp.~553--572.

\bibitem{polik2007survey}
{\sc I.~P{\'o}lik and T.~Terlaky}, {\em A survey of the {{S}}-lemma}, SIAM
  Review, 49 (2007), pp.~371--418.

\bibitem{pong2014generalized}
{\sc T.~K. Pong and H.~Wolkowicz}, {\em The generalized trust region
  subproblem}, Computational Optimization and Applications, 58 (2014),
  pp.~273--322.

\bibitem{rendl1997semidefinite}
{\sc F.~Rendl and H.~Wolkowicz}, {\em A semidefinite framework for trust region
  subproblems with applications to large scale minimization}, Mathematical
  Programming, 77 (1997), pp.~273--299.

\bibitem{salahi2018efficient}
{\sc M.~Salahi and A.~Taati}, {\em An efficient algorithm for solving the
  generalized trust region subproblem}, Computational and Applied Mathematics,
  37 (2018), pp.~395--413.

\bibitem{stern1995indefinite}
{\sc R.~J. Stern and H.~Wolkowicz}, {\em Indefinite trust region subproblems
  and nonsymmetric eigenvalue perturbations}, SIAM Journal on Optimization, 5
  (1995), pp.~286--313.

\bibitem{sturm2003cones}
{\sc J.~F. Sturm and S.~Zhang}, {\em On cones of nonnegative quadratic
  functions}, Mathematics of Operations Research, 28 (2003), pp.~246--267.

\bibitem{wang2016linear}
{\sc J.~Wang and Y.~Xia}, {\em A linear-time algorithm for the trust region
  subproblem based on hidden convexity}, Optimization Letters,  (2016),
  pp.~1--8.

\bibitem{yakubovich1971sprocedure}
{\sc V.~A. Yakubovich}, {\em S-procedure in nonlinear control theory}, Vestnik
  Leningrad University, 1 (1971), pp.~62--77.

\bibitem{yang2013two}
{\sc B.~Yang and S.~Burer}, {\em A {{Two}}-{{Variable Approach}} to the
  {{Two}}-{{Trust}}-{{Region Subproblem}}}, tech. report, Tech. Report,
  Department of Mathematics, University of Iowa, 2013.

\bibitem{ye1992new}
{\sc Y.~Ye}, {\em A new complexity result on minimization of a quadratic
  function with a sphere constraint}, in Recent {{Advances}} in {{Global
  Optimization}}, {C. Floudas and P. Pardalos}, ed., {Princeton University
  Press}, 1992, pp.~19--31.

\bibitem{ye2003new}
{\sc Y.~Ye and S.~Zhang}, {\em New results on quadratic minimization}, SIAM
  Journal on Optimization, 14 (2003), pp.~245--267.

\bibitem{zhang2010derivative}
{\sc H.~Zhang, A.~R. Conn, and K.~Scheinberg}, {\em A derivative-free algorithm
  for least-squares minimization}, SIAM Journal on Optimization, 20 (2010),
  pp.~3555--3576.

\end{thebibliography}

\end{document}